\documentclass[11pt]{article}
\usepackage[a4paper]{anysize}\marginsize{3.5cm}{3.5cm}{1.3cm}{2cm}
\pdfpagewidth=\paperwidth \pdfpageheight=\paperheight
\usepackage{amsfonts,amssymb,amsthm,amsmath,eucal}
\usepackage{pgf}
\usepackage{bbm,array}
\usepackage{pstricks}
\usepackage{tikz} 

\theoremstyle{plain}
\newtheorem{thm}{Theorem}[section]
\newtheorem{theorem}[thm]{Theorem}

\newtheorem{lemma}[thm]{Lemma}
\newtheorem{proposition}[thm]{Proposition}
\newtheorem{corollary}[thm]{Corollary}

\theoremstyle{definition}
\newtheorem{definition}[thm]{Definition}

\newtheorem{example}[thm]{Example}

\newtheorem{thevarthm}[thm]{\varthmname}

\newenvironment{varthm*}[1]{\trivlist\item[]{\bf #1.}\it}{\endtrivlist}

\renewcommand\geq{\geqslant}

\renewcommand\leq{\leqslant}

\newcommand\be{\begin{eqnarray*}}
\newcommand\ee{\end{eqnarray*}}

\newcommand\F{\mathbb F}
\renewcommand\P{\mathbb P}

\newcommand\calo{{\mathcal O}}
\newcommand\cali{{\mathcal I}}

\newcommand\newop[2]{\def#1{\mathop{\rm #2}\nolimits}}
\newop\log{log}
\newop\ord{ord}
\newop\Gal{Gal}
\newop\SL{SL}
\newop\Bl{Bl}
\newop\mult{mult}
\newop\mass{mass}
\newop\div{div}
\newop\codim{codim}
\newop\sing{sing}
\newop\vdim{vdim}
\newop\edim{edim}
\newop\Ass{Ass}
\newop\size{size}
\newop\reg{reg}
\newop\satdeg{satdeg}
\newop\supp{supp}
\newcommand\eqnref[1]{(\ref{#1})}

\newcommand\alphahat{\widehat\alpha}

\def\keywordname{{\bfseries Keywords}}%
\def\keywords#1{\par\addvspace\medskipamount{\rightskip=0pt plus1cm
\def\and{\ifhmode\unskip\nobreak\fi\ $\cdot$
}\noindent\keywordname\enspace\ignorespaces#1\par}}
\def\subclassname{{\bfseries Mathematics Subject Classification
(2000)}\enspace}
\def\subclass#1{\par\addvspace\medskipamount{\rightskip=0pt plus1cm
\def\and{\ifhmode\unskip\nobreak\fi\ $\cdot$
}\noindent\subclassname\ignorespaces#1\par}}

\begin{document}

\author{S.~Di Rocco\footnote{Di Rocco's research was partially support by the VR grant NT:2010-5563.} , A.~Lundman, T.~Szemberg\footnote{Szemberg research was partially supported
by NCN grant UMO-2011/01/B/ST1/04875}}
\title{The effect of points fattening on Hirzebruch surfaces}
\date{\today}
\maketitle
\thispagestyle{empty}

\begin{abstract}
   The purpose of this note is to study initial
   sequences of $0$--dimensional subschemes of
   Hirzebruch surfaces and classify subschemes
   whose initial sequence has the minimal possible growth.
\subclass{MSC 14C20 \and MSC 13C05 \and MSC 14N05 \and MSC 14H20 \and MSC 14A05}
\end{abstract}

\section{Introduction}
\label{sec:intro}
   Zeroschemes of fat points on algebraic varieties consist of a finite number of points given
   with prescribed multiplicity. They play a fundamental role in interpolation problems for algebraic curves
   and in the theory of minimal resolutions of ideals, see for example \cite{Cat91}, \cite{CCMO03}, \cite{CHT11}.

    Fat points subschemes on the projective plane  are best  (but by far not completely) understood.
   In \cite{BocCha11} Bocci and Chiantini introduced a new invariant called the    \emph{initial degree}.
     For a $0$--dimensional reduced subscheme $Z$ of $\P^2$  the initial degree $\alpha(Z)$ of $Z$ (or rather of the homogeneous ideal $I_Z$
   defining $Z$) is  the minimal degree of a non--zero
   element in $I_Z$.  In other words $\alpha(Z)$ is the minimal
   integer $d$ such that there exists a plane curve of degree $d$
   vanishing at all points of $Z$. Similarly, for a positive integer $m$
   one defines $\alpha(mZ)$ to be the minimal integer $d$ such that there
   exists a plane curve of degree $d$ vanishing to the order at least $m$ in all
   points of $Z$. Computing $\alpha(mZ)$ is in general a challenging problem, governed partially by the celebrated Nagata Conjecture \cite{Nag59}.
   Bocci and Chiantini studied subschemes $Z$ of $\P^2$ with
   $$\alpha(2Z)=\alpha(Z)+1$$
   and obtained a full classification \cite[Theorem 1.1]{BocCha11}.
   This result has already motivated a considerable amount of research,
   see e.g. \cite{BauSze13}, \cite{Lanckorona}, \cite{DST13}, \cite{Jan13}.
   Whereas \cite{DST13} deals still with fat points in the projective
   plane, \cite{BauSze13} contains results paralleling those of
   \cite{BocCha11} for $\P^3$ and formulates a conjectural picture
   for higher dimensional projective spaces and higher dimensional
   linear subschemes of $\P^n$. That approach has been taken on
   by Jannsen \cite{Jan13}, who studies initial degrees for ACM
   unions of lines in $\P^3$. The article \cite{Lanckorona} extends
   Bocci and Chiantini results from the homogeneous to the multi-homogeneous
   case, more precisely it deals with $0$--dimensional subschemes
   of $\P^1\times\P^1$. Necessary modifications of the original
   definition of $\alpha(Z)$ taken on in \cite{Lanckorona} suggest
   that the fattening effect problem is not restricted to projective
   spaces or their products. In particular that article suggests the
   following definition.
\begin{definition}[Initial degree and initial sequence]\label{def:initial}
   Let $X$ be a smooth projective variety with an ample line bundle $L$
   and let $Z$ be a reduced subscheme of $X$ defined by the ideal sheaf $\cali_Z\subset\calo_X$.
   For a positive integer $m$
   we define the \emph{initial degree (with respect to $L$)} of the subscheme $mZ$ as
   $$\alpha(mZ):=\min\left\{d:\; H^0(X;dL\otimes\cali_Z^m)\neq 0\;\right\}.$$
   The \emph{initial sequence (with respect to $L$)} of $Z$ is then the sequence
   $$\alpha(Z), \alpha(2Z), \alpha(3Z),\ldots$$
\end{definition}
   This definition leaves some ambiguity in choosing a polarization $L$.
   If $X$ is the projective space, then one works naturally with the
   ample generator $\calo_{\P^n}(1)$, see \cite{BocCha11} and \cite{BauSze13}.
   Similarly for products of projective spaces the natural polarisation is ${\mathcal O}_{\P^1\times\P^1}(1,1) .$ In both cases these polarisations
   can be considered the ``minimal ones" as  their
   degree and volume are minimal.
   We believe that working with ample classes which are minimal with respect to their
   degree and volume  is the most natural and effective way  to study points fattening effect
  on any class of algebraic varieties. While this does not lead to a unique choice
   for certain varieties, e.g. abelian varieties may carry a lot of principal
   polarizations \cite{Lan87}, there are interesting classes of
   varieties where the choice is unique, e.g. projective bundles over $\P^1$.

   In this paper we take the first step towards
   understanding how the initial sequence determines
   the geometry of zero dimensional subschemes of Hirzebruch surfaces.
   Our results may be summarized as follows, see the next section for precise
   explanation of the notation.
\begin{varthm*}{Theorem}
   Let $Z\subset \F_r$ be a finite set of reduced points and let $L=L_r=(r+1)F_r+E_r$.
   \begin{itemize}
      \item[a)] The initial sequence
      $$\alpha(Z)=\dots=\alpha((r+1)Z)=\alpha((r+2)Z)=\alpha((r+3)Z)=a$$
      is not admissible.
      \item[b)] If $Z$ has the initial sequence
      $$\alpha(Z)=\dots=\alpha((r+1)Z)=\alpha((r+2)Z)=a$$
    then $Z$ is a single point contained in the negative curve $E_r$
      and moreover $a=1$.
      \item[c)] If $Z$ has the initial sequence
      $$\alpha(Z)=\dots=\alpha((r+1)Z)=a$$
     then:
      \subitem either all points in $Z$ are contained in a single fiber of the projection $\varphi_r:\F_r\to\P^1$ and $a=1$;
      \subitem or $r=1$ and  the points in $Z$ are all intersection points of a general star configuration of $a$ lines
      and $a$ lines passing through one point. In this situation all positive values of $a$ can occur.
   \end{itemize}
\end{varthm*}
   This result is proved in Propositions \ref{prop:no r+2 jumps by 0} and \ref{prop:r+1 jumps by 0}
   and Theorem \ref{thm:r jumps by 0}.

   We conclude the introduction with the definition of an asymptotic counterpart of the initial degree.
\begin{definition}[Waldschmidt constant]\label{def:waldschmidt}
   Keeping the notation from Definition \ref{def:initial} we define
   the \emph{Waldschmidt constant} of $Z$ (\emph{with respect to $L$}) as
   the limit
   $$\alphahat(Z):=\lim_{m\to\infty}\frac{\alpha(mZ)}{m}.$$
\end{definition}
   The existence of the limit in the definition follows
   from the subadditivity of the initial sequence by a standard argument based on the
   Fekete Lemma \cite{Fek23}. It follows also
   that $\alphahat(Z)$ is in fact the infimum of the sequence terms.

   Waldschmidt constants are interesting invariants that
   were recently rediscovered and studied by Harbourne,
   see e.g. \cite{BocHar10b}. We show in
   Proposition \ref{thm:chudnovsky} that somewhat unexpectedly
   they behave very regularly on Hirzebruch surfaces.

\section{Hirzebruch surfaces}\label{sec:hirz surf}
   In this section we collect some general facts about Hirzebruch surfaces.
   Recall that for a non--negative integer $r$, the Hirzebruch surface
   $\F_r$ is defined as the projectivization of the vector
   bundle $\calo_{\P^1}\oplus\calo_{\P^1}(r)$. The effective
   cone of $\F_r$ is spanned by classes of two smooth rational
   curves: $F_r$, a fiber of the projection $\varphi_r: \F_r\to\P^1$
   and $E_r$, a section of $\varphi_r$ satisfying $E_r^2=-r$.
   If $r=0$, then $\F_0=\P^1\times\P^1$. The effect of points fattening
   on $F_0$ was studied in \cite{Lanckorona}, hence here we focus our attention on the case $r\geq 1$.
   In this situation $E_r$ is the unique curve $C$ on $\F_r$ with $C^2=-r$
   and this is also the unique curve $C$ on $\F_r$ with $h^0(\F_r,C)=1$.

   In this paper we study the effect of fattening with respect to the line bundle $L_r=(r+1)F_r+E_r$ on $\F_r$,
   which is the ample line bundle of minimal
   degree (self--intersection) and volume. It is worth to note
   that also from the point of view of toric geometry, working with
   $L_r$ is a natural choice. Indeed, $L_r$ is a line bundle
   associated to the minimal polytope whose normal fan is
   isomorphic to the normal fan of $\F_r$. By minimality we mean here
   minimal euclidian volume of the polytope or equivalently
   minimal edge-length. Note that all these properties characterise also the minimal class $\calo_{\P^2}(1)$ on $\P^2$.

   For $r=1$ the surface $\F_1$ is just a blowing up
   $\pi_1:\F_1\to\P^2$ of the projective plane
   $\P^2$ in a point $Q$. The section $E_1$ is the
   exceptional divisor of this blowing up.

   The surfaces $\F_r$ are related by birational transformations called \emph{elementary transformations}.
 Blowing up a point  $P\in\F_r$ on $ E_r$ and then contracting the preimage of the fiber through that point gives a brational map to $\F_{r+1}.$ If we blow up instead a point not on  on $ E_r$ and then contract the preimage of the fiber through that point we obtain a  brational map to $\F_{r-1}.$
\section{An auxiliary Lemma}\label{sec:general}
   In this section we present a lemma concerning properties
   of singularities of plane curves.
\begin{lemma}\label{lem:max number of singular points}
   Let $C$ be a reduced divisor in $\P^2$ of degree $d$ with a point $Q$
   of multiplicity $m\geq 2$. Then $C$ has at most
   $\binom{d}{2}-\binom{m}{2}$ singular points away of $Q$.

   Moreover if the number of singular points is exactly equal to
   $\binom{d}{2}-\binom{m}{2}$, then $C$ is a union of
   $m$ lines passing through $Q$ and a set of $d-m$
   general lines with $\binom{d-m}{2}$ intersection points,
   none of them lying on a line from the first set of lines (i.e. those passing through $Q$).
\end{lemma}
\proof
   We proceed by induction on the multiplicity $m$. For $m=2$
   and $d$ arbitrary, the assertion follows from \cite[Lemma 7.5]{GMS06}.

   We assume now that the Lemma holds for $m$ and arbitrary $d$
   and we want to conclude that it holds for $m+1$ and arbitrary $d$.
   Let $C$ be a divisor of degree $d$ passing through the point $Q$
   with multiplicity $m+1$. We consider
   two cases.

   \textbf{Case 1.} We assume that there is only one irreducible component $\Gamma$ of $C$
   of degree $\gamma\leq d$ passing through the point $Q$. In particular,
   $\Gamma$ is not a line since we have $m\geq 2$.
   By adjunction we have then that the number $s_{\Gamma}$
   of singular points on $\Gamma$ is at most
   $$s_{\Gamma}\leq \frac{(\gamma-1)(\gamma-2)}{2}-\binom{m+1}{2}.$$
   By induction assumption, there are at most
   $\binom{d-\gamma}{2}$ singular points on the residual curve $C-\Gamma$.
   There are at most $\gamma\cdot(d-\gamma)$ additional singular points
   coming from intersection points of $\Gamma$ and $C-\Gamma$. Altogether
   the number $s_C$ of singular points on the divisor $C$ is in this case bounded by
   $$s_C\leq \frac{(\gamma-1)(\gamma-2)}{2}-\binom{m+1}{2}+\binom{d-\gamma}{2}+\gamma\cdot(d-\gamma).$$
   This number is strictly less than $\binom{d}{2}-\binom{m+1}{2}$ since we can assume that $\gamma\geq 3$.

   \textbf{Case 2.} Now assume that there are  at least two  divisors
   $\Gamma_1$ and $\Gamma_2$ passing through the point $Q.$ 
   Note that we do not assume here that $\Gamma_i$ are irreducible.
   The point is that we have 
   $$1\leq m_i=\mult_Q(\Gamma_i)\leq m,$$
   so that we can apply the induction assumption to both divisors $\Gamma_1$ and $\Gamma_2$. 
   Then by
   this assumption and the same counting as in Case 1. we obtain
   \begin{equation}\label{eq:bounding s_C}
   s_C\leq \binom{\gamma_1}{2}-\binom{m_1}{2}+\binom{\gamma_2}{2}-\binom{m_2}{2}
   +(\gamma_1\cdot\gamma_2-m_1\cdot m_2),
   \end{equation}
   where $\gamma_i$ is the degree of $\Gamma_i$ (so that $\gamma_1+\gamma_2=d$). Hence
   $$s_C\leq \binom{\gamma_1+\gamma_2}{2}-\binom{m_1+m_2}{2}= \binom{d}{2}-\binom{m+1}{2},$$
   which proves the first assertion. 

   Notice that the inequality in \eqnref{eq:bounding s_C} is strict
   unless all summands are maximal, i.e.
   $$s_{\Gamma_i}=\binom{\gamma_1}{2}-\binom{m_i}{2}$$
   and $\Gamma_1$ intersects $\Gamma_2$ away of $Q$ in exactly $\gamma_1\gamma_2-m_1m_2$ points.
   By induction this implies that $\Gamma_1$ and $\Gamma_2$ split into lines and the configuration
   of these lines satisfies the last assertion of the Lemma.
\endproof

\section{Subschemes with initial sequences of minimal growth}\label{sec:minimal}
   Sections of $\calo_{\P^2}(d)$ are just homogeneous polynomials of degree $d$.
   If such a section vanishes at a point $P\in\P^2$ to the order $m\geq 2$,
   then any of its directional derivatives vanishes at this point to the order
   $\geq m-1$. This shows that on $\P^2$ the relation $\alpha((m+1)Z)>\alpha(mZ)$ always holds.

   On the other hand on $\F_0=\P^1\times\P^1$, it might easily happen that
   $\alpha((m+1)Z)=\alpha(mZ)$ for some subscheme $Z\subset\F_0$. However
   two consecutive equalities of this kind are not possible,
   see \cite[Corollary 2.4]{Lanckorona}.

   The reason behind these two statements is of the same flavor.
   In the case of $\P^2$ there are no singular divisors in the linear
   system $|\calo_{\P^2}(1)|$ (because they are all lines),
   whereas there are divisors with multiplicity at most $2$
   in the minimal polarization $L_0$. This pattern extends to all Hirzebruch surfaces.
\begin{proposition}\label{prop:no r+2 jumps by 0}
   There is no zero--dimensional subscheme $Z\subset \F_r$ with
   vanishing sequence
   \begin{equation}\label{eq:r+2 jums by 0}
      \alpha(Z)=\alpha(2Z)=\ldots=\alpha((r+3)Z).
   \end{equation}
\end{proposition}
\proof
   Let $a:=\alpha(mZ)$. If $a=1$, then it is elementary to check
   (for example intersecting with the fiber through the given point)
   that there is no section in $H^0(\F_r, L_r)$ vanishing at a point
   to order $\geq r+3$. Hence it must be $\alpha((r+3)Z)\geq 2$
   contradicting \eqnref{eq:r+2 jums by 0}.

   If $a\geq 2$, then the assertion follows from inequality (\ref{eq:chudnovsky all cases})
   below.
\endproof
   The following Example shows that the statement in Proposition \ref{prop:no r+2 jumps by 0}
   is optimal.
\begin{example}\label{ex:point}
   Let $Z$ be subscheme of $\F_r$ consisting of a single point $P\in E_r$.
   Then
   \begin{equation}\label{eq:r+1 jumps by 0}
      \alpha(Z)=\alpha(2Z)=\ldots=\alpha((r+2)Z)=1.
   \end{equation}
\end{example}
\proof
   Let $F_P$ denote the fiber passing through the point $P$.
   Then the divisor $(r+1)F_P+E_r\in |L_r|$
   has multiplicity $r+2$ at $P$. One easily  checks
   that there are no divisors in $|L_r|$ with higher multiplicity at $P$.
\endproof
   We will show that the situation in Example \ref{ex:point}
   is the unique situation in which \eqnref{eq:r+1 jumps by 0} holds.
   To this end we prove first the following generalization of a result of Chudnovsky
   \cite{Chu81}. It puts some universal constrains on the order of
   growth of the initial sequence $\alpha(mZ)$. This kind of results
   for ideals of points in the projective plane was first obtained by
   Skoda and Waldschmidt. Our approach here is modeled on the algebraic
   proof for $\P^2$ by Harbourne and Huneke \cite[Proposition 3.1]{HaHu13}.
   This result might be of independent interest. As usually, we assume $r\geq 1$.
\begin{theorem}\label{thm:chudnovsky}
   Let $Z\subset\F_r$ be a reduced $0$--dimensional subscheme.
   Then the inequality
   \begin{equation}\label{eq:chudnovsky all cases}
      \frac{\alpha(mZ)}{m}\geq \frac{\alpha(Z)}{r+2}
   \end{equation}
   holds for all $m\geq 1$.
\end{theorem}
\proof
   If $\alpha(Z)=1$, then $Z$ is contained in a divisor
   $D\in|L_r|$. As in Example \ref{ex:point} one can show
   that then $\alpha(k(r+2)Z)\geq k$.
   Indeed, this follows from the fact that
   there is no divisor with a point of multiplicity
   $m\geq k(r+2)+1$ in $kL_r$ for all $k\geq 1$. Hence
   \eqnref{eq:chudnovsky all cases} holds in this case.

   Thus we assume that $a:=\alpha(Z)\geq 2$. There is a subset $W\subset Z$ of
   points $\left\{Q_1,\ldots,Q_t\right\}$ with minimal  $t$ such that
   $\alpha(W)=a$ (of course it might happen that $W=Z$). Since
   vanishing at a point is a single linear condition on sections
   of a linear series, it must be in fact
   $$t=\dim(H^0(\F_r,(a-1)L_r))=\frac12ra(a-1)+a^2.$$
   \textbf{Claim.} The linear system $|aL_r\otimes\cali_W|$
   has either no fixed components or it has exactly one fixed component based
   on the negative curve $E_r$ and in this case the set $W$ is disjoint from $E_r$.

   In order to prove the claim we assume that $\Gamma$ is a base component
   of $|aL_r\otimes\cali_W|$. By the choice of $W$, for every point $Q\in W$
   there exists a divisor $C_Q\in |(a-1)L_r\otimes \cali_{W\setminus\left\{Q\right\}}|$
   not vanishing at $Q$ and vanishing at all other points in $W$.

   Now, if $\Gamma$ is not a component of $C_Q$, then let $C'\in|L_r\otimes\cali_Q|$
   be a divisor vanishing at $Q$ and such that $\Gamma$ is not its component
   (such a divisor can be found since $L_r$ is very ample). Then
   $C'+C_Q\in|aL_r\otimes\cali_W|,$
   is a divisor not containing $\Gamma$, a contradiction.

   Hence $\Gamma$ is a component of all divisors $C_Q$. This implies $\Gamma\cap W=\emptyset$.
   Indeed, if there were a point $R\in\Gamma\cap W$, then we would have $R\in\Gamma\subset C_R$,
  contradicting  the definition of $C_R$.

   If $\dim(H^0(\Gamma))\geq 2$, then there is a divisor $\Gamma'$ linearly equivalent to $\Gamma$
   passing through a given point $Q\in W$. Then
   $$C_Q-\Gamma+\Gamma'\in |(a-1)L_r\otimes\cali_W|$$
   is an effective divisor (since $\Gamma$ is a component of $C_Q$)
   vanishing at all points of $W$. This contradicts $\alpha(W)=a$.

   If $\dim(H^0(\Gamma))=1$, then $\Gamma=bE_r$ for some $b\geq 1$.
   Indeed, $\dim(H^0(bE_r+cF))\geq c+1$ for all $b,c\geq 0$
   (it suffices to count reducible divisors in that linear series).
   Thus the Claim above is established.

   Now, let $A$ be a divisor in the linear system $|\alpha(mZ)\cdot L_r|$
   vanishing along $Z$ with multiplicities $\geq m$.

   By the Claim, there exists a divisor $B$ in $aL_r\otimes\cali_W$ either with no common
   component with $A$ or with a common component supported on $E_r$ in which case
   no point from $W$ lies on that component. Hence an obvious modification of Bezout
   Theorem implies that
   $$a\alpha(mZ)L_r^2=A\cdot B\geq \left(\frac12ra(a-1)+a^2\right)m$$
   It follows that
   \begin{equation}\label{eq:chudnovsky with r}
      \frac{\alpha(mZ)}{m}\geq \frac{\frac12r(a-1)+a}{r+2},
   \end{equation}
   which implies
   \begin{equation}\label{eq:chudnovsky proof}
      \frac{\alpha(mZ)}{m}\geq \frac{a}{r+2}
   \end{equation}
   and completes the proof of the Theorem.
\endproof
   From the proof of Theorem \ref{thm:chudnovsky} we get immediately the following
   useful corollary.
\begin{corollary}\label{cor:chudnovsky}
   Let $Z\subset \F_r$ be a reduced $0$--dimensional subscheme with $\alpha(Z)\geq 2$.
   Then
   \begin{equation}\label{eq:chudnovsky alpha geq 2}
      \frac{\alpha(mZ)}{m}> \frac{\alpha(Z)}{r+2}
   \end{equation}
   holds for all $m\geq 1$.
\end{corollary}
   Now we are in the position to prove the equivalence in Example \ref{ex:point}.
\begin{proposition}\label{prop:r+1 jumps by 0}
   Let $Z\subset\F_r$ be a $0$--dimensional subscheme with
   \begin{equation}\label{eq:r+1 jumps by 0 again}
      \alpha(Z)=\alpha(2Z)=\ldots=\alpha((r+2)Z).
   \end{equation}
   Then $Z$ is a single point $P\in E_r$ and $\alpha(Z)=1$.
\end{proposition}
\proof
   It follows immediately from \eqnref{eq:chudnovsky alpha geq 2} that
   \eqnref{eq:r+1 jumps by 0 again} implies $\alpha(Z)=1$. Then it is
   elementary to check that any divisor in $|L_r|$ has at most one
   point with multiplicity $\geq r+2$ and this point must then lie
   on $E_r$.
\endproof
   In the rest of the paper we study subschemes $Z$ for
   which the first $r+1$ terms in the initial sequence are constant.
   Our main result in this
   case is the following theorem.
\begin{theorem}\label{thm:r jumps by 0}
   Let $Z\subset \F_r$ be a $0$--dimensional subscheme.
   If
   $$
   \alpha(Z)=\alpha(2Z)=\ldots=\alpha((r+1)Z))=:\alpha
   $$
   then\\
   a) $Z$ is contained in a single fiber $F$ if $r\geq 2$;\\
   b) for $r=1$ the set $Z$ is the preimage under the blow-up map $\pi_1:\F_1\mapsto\P^2$
       of the intersection points of $\alpha$
      lines passing through the point $Q=\pi_1(E_1)$ and a star configuration
      of additional $\alpha$ general lines such that no intersection
      point of the star configuration lies on a line passing through $Q$.\\
   The figure below shows a configuration described in part b)
   for $\alpha=3.$
\end{theorem}

\begin{figure}[h!]
\centering
\begin{tikzpicture}
\draw  [thick] (0,0)--(5,0);
\draw  [thick] (0,1)--(5,-2);
\draw  [thick] (0,-1)--(5,2);
\draw  [thick] (2,-3)--(4,4);
\draw  [thick] (1,2)--(5,-3);
\draw  [thick] (3.5,4)--(5,-3.5);

\draw [thick, red] (5/3,0) circle (0.01) node[below]{$Q$};
\draw [thick, red] (13/5,0) circle (0.09);
\draw [thick, red] (85/37,14/37) circle (0.09);
\draw [thick, red] (45/13,-14/13) circle (0.09);
\draw [thick, red] (20/7,0) circle (0.09);
\draw [thick, red] (110/41,-25/41) circle (0.09);
\draw [thick, red] (90/29,25/29) circle (0.09);
\draw [thick, red] (43/10,0) circle (0.09);
\draw [thick, red] (225/56,79/56) circle (0.09);
\draw [thick, red] (205/44,-79/44) circle (0.09);

\draw [thick] (53/19,-9/38) circle (0.09);
\draw [thick] (73/15,-17/6) circle (0.09);
\draw [thick] (63/17,101/34) circle (0.09);
\end{tikzpicture}
\caption{$a$ lines in a pencil and a star configuration of $a$ lines}
\end{figure}
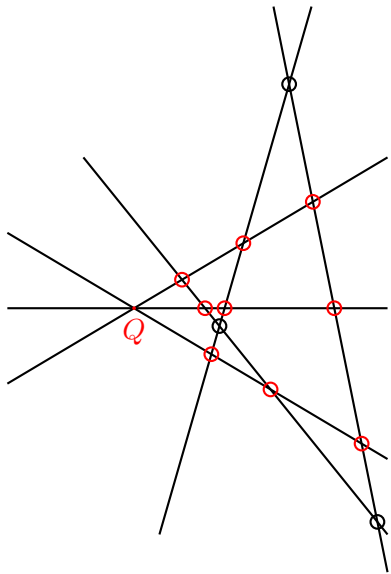


\proof
   Note  that the inequality \eqref{eq:chudnovsky with r} with $m=r+1$
   gives
   $$\alpha\leq\frac{r(r+1)}{r^2+r-2},$$
   after rearranging terms. For $r\geq 2$ this implies $\alpha=1$.

   Now, for $r\geq 2$,
   let $D\in |L_r|$ be a divisor with $\mult_PD\geq r+1$ for all points $P\in Z$.
   Let $F$ be the fiber of the projection $\varphi_r:\F_r\to P_1$ passing
   through a point $P\in Z$. Since $D\cdot F=1$ and $P$ is a point of multiplicity
   at least $r+1$ on $D$, it follows from Bezout Theorem that $F$ is at least
   an $r$--fold component of $D$. This argument holds for any other point $Q\in Z$.
   Since $D-2F_r$ is obviously not an effective divisor, it follows that $P$
   and $Q$ (and hence any other point in $Z$) must lie in the same fiber.

   We conclude the proof with an argument for $r=1$.
   If $\alpha=1$, then we are in the previous case, so we assume
   $\alpha\geq 2$.
   It is convenient to change the notation a little bit.
   We denote as usual $H=\pi_1^*(\calo_{\P^2}(1))$ the
   pullback of the hyperplane bundle and by $E$ the exceptional
   divisor of the blow up $\pi_1:\F_1\to\P^2$. Let $Q=\pi_1(E_1).$   Let $D\in |\alpha L_1|=|\alpha(2F_1-E_1)|$ be
   a divisor with $\mult_PD\geq r+1$ for all points $P\in Z$.
   The existence of $D$ implies that there exists a plane divisor $C$
   of degree $d=2\alpha$ with a point of multiplicity at least $m=\alpha$ at $Q$
   and multiplicity at least $\alpha$ in all other points in $Z$.
   The assumption $\alpha(Z)=\alpha$ implies in turn that there is no plane curve
   of degree $2(\alpha-1)$ with multiplicity at least $\alpha-1$ at $Q$ passing through
   all other points in $Z$. A naive dimension count shows that
   $$\binom{2\alpha}{2}-\binom{\alpha}{2}-{\rm supp}(Z)\leq 0,$$
   we obtain then that there are at least $\binom{2\alpha}{2}-\binom{\alpha}{2}$ points in $Z$.
   Then Lemma \ref{lem:max number of singular points} implies the
   assertion of the Theorem.
\endproof

\paragraph*{\emph{Acknowledgement.}}
   This paper has been written while the third named author visited KTH.
   It's a pleasure to thank the department of Mathematics of KTH for
   providing excellent working conditions. We are grateful to the G\"oran Gustafsson Foundation for financial support.



\bigskip \small

\bigskip
   Sandra Di Rocco,
   Department of Mathematics, KTH, 100 44 Stockholm, Sweden.

\nopagebreak
   \textit{E-mail address:} \texttt{dirocco@math.kth.se}

\bigskip
   Anders Lundman,
   Department of Mathematics, KTH, 100 44 Stockholm, Sweden.

\nopagebreak
   \textit{E-mail address:} \texttt{alundman@kth.se}

\bigskip
   Tomasz Szemberg,
   Instytut Matematyki UP,
   Podchor\c a\.zych 2,
   PL-30-084 Krak\'ow, Poland.

\nopagebreak
   \textit{E-mail address:} \texttt{tomasz.szemberg@gmail.com}


\end{document}